\documentclass[12pt]{amsart}
\usepackage[utf8]{inputenc}
\usepackage[T2A,T1]{fontenc}
\usepackage[russian,english]{babel}
\usepackage{textcase}

\usepackage[all, cmtip]{xy}
\usepackage{amsfonts,amsmath,amstext,amsbsy,amssymb, amsopn,amsthm,upref,eucal, mathrsfs}

\newtheorem{proposition}{Proposition}
\newtheorem{theorem}{Theorem}
\newtheorem{lemma}{Lemma}
\newtheorem{corollary}{Corollary}
\newtheorem{remark}{Remark}

\begin{document}

\title[Decay of Correlations for Multiplicative Functionals]{Decay of Correlations for Multiplicative Functionals of Determinantal Point Processes}

\author{Alexander\ I.~Bufetov}
\address{Steklov Mathematical Institute of Russian Academy of Sciences, Moscow}
\email{bufetov@mi-ras.ru}

\thanks{This work was carried out at the International Mathematical Center of the Steklov Mathematical Institute (MIAN) with financial support from the Ministry of Science and Higher Education of the Russian Federation (Agreement No. 075-15-2025-303).}

\date{}

\begin{abstract}
The main result of this paper is an upper bound on the correlations of regularized multiplicative functionals of determinantal point processes. In particular, for the sine process, the correlations of the multiplicative functional associated with a rapidly decaying function decay quadratically. It follows that the spectral density of such functionals belongs to the Sobolev space $H_1$.
\end{abstract}

\maketitle

\section{Introduction and statement of the main result}
\subsection{Multiplicative functionals of determinantal \\ point processes}
Let $(E,d)$ be a complete separable metric space,
let $\mu$ be a $\sigma$-finite measure on $E$, and let
$\varPi: E\times E \to \mathbb{C}$ be a Borel function of two variables that is bounded on compact sets,
and such that the operator $\varPi f(x) = \int\varPi(x,y)f(y)d\mu(y)$,
initially defined on bounded Borel functions with compact support, extends to
an orthogonal projection acting on $L_2(E,\mu)$ (all Hilbert spaces considered in this paper are complex).
By the Macchi–Soshnikov–Shirai–Takahashi theorem \cite{Macchi,Soshnikov-2,Shirai-Tokohashi},
the projection $\varPi$ determines a determinantal measure on the space $\operatorname{Conf}(E)$ of locally finite point configurations on $E$.

If $f$ is a bounded Borel function with compact support on $E$,
and $\Psi_{1+f}$ is the multiplicative functional associated with $f$ on the space of configurations $\operatorname{Conf}(E)$ defined by the formula
\[
\Psi_{1+f}(X) = \prod_{x \in X}(1+f(x)),
\]
then (cf. \cite{Bufetov-Cond,Bufetov-Symm}) the determinantal measure $\mathbb P_{\varPi}$ is uniquely determined by the condition
\[
\mathbb{E}_{\mathbb P_{\varPi}} \Psi_{1+f} = \det(1+\varPi f\varPi)= \det(1+\sqrt{f}\varPi\sqrt{f}).
\]
Examples of measures that would subsequently be called determinantal already appear in Andr\'ei\'ef's work \cite{Andreief}.

The notation $\sqrt{f}\varPi\sqrt{f}$ is well-defined for both real- and complex-valued functions $f$. Either square root of $f$ may be chosen; the only requirement is that the same square root be used on the left and on the right hand side. Now set
\[
\widetilde{\Psi}_{1+f} = \Psi_{1+f} \cdot \exp\bigl(-\mathrm{tr}(\varPi f\varPi)\bigr) = \Psi_{1+f}\exp(-\mathrm{tr}\sqrt{f}\varPi\sqrt{f}).
\]

We next introduce the normalized multiplicative functionals. Let $HS$ to be the space of Hilbert–Schmidt operators
acting on $L_2(E,\mu)$. Consider the space $HS(\varPi)$ of Borel functions $f$ for which the operator $f\varPi$ is Hilbert–Schmidt, and equip this space with the Euclidean structure
defined by $\|f\|^2_{HS(\varPi)} = \|f\varPi\|^2_{HS}$; by definition, we have
\begin{equation*}
\langle f, h \rangle_{HS(\varPi)} = \mathrm{tr}(f\varPi \overline{ h}) = \langle f\varPi, h\varPi \rangle_{HS}.
\end{equation*}

For the Hilbert space obtained by completing the space of bounded Borel functions with compact support with respect to the Euclidean structure ${<,>_{HS(\varPi)}}$, we retain, with a slight abuse of notation, the symbol $HS(\varPi)$.

Recall that the regularized Hilbert–Carleman determinant
\cite{carleman} is defined for a finite-rank bounded operator $A$ by the formula
$$\det\nolimits_2(1+A)=\det(1+A)\exp(-\mathrm{tr} A)$$
and is then extended by continuity to the entire space of Hilbert–Schmidt operators.
\begin{proposition}
The map $f \to \widetilde{\Psi}_{1+f}$ extends to a \\ continuous map
$HS(\varPi)\to L_1(\operatorname{Conf}(E), \mathbb P_{\varPi})$, and for every $f \in HS(\varPi)$,
we have:
\[
\mathbb{E}_{\mathbb P_{\varPi}} \widetilde{\Psi}_{1+f} = \det\nolimits_2(1+\varPi f\varPi)=\det\nolimits_2(1+f\varPi).
\]
\end{proposition}
\subsection{A general theorem on decay of correlation}
The main result of this paper is an upper bound on the covariance of two normalized multiplicative functionals.

We use $L_2(\Pi)$ to denote the space of square-integrable functions with respect to the measure $\Pi(x,x) d\mu(x)$.

We will also need another embedding of the space of bounded Borel
functions with compact support into the space of Hilbert–Schmidt operators,
namely $f\to\varPi f\varPi$. Define the Hilbert space
$H^{(2)}(\varPi)$ as the completion of the space of bounded
Borel functions with compact support in the norm
\begin{equation*}
\|f\|^2_{H^{(2)}(\varPi)}=\|\varPi f\varPi\|^2_{HS};
\end{equation*}
defined by the inner product
\begin{equation*}
\langle f,h\rangle_{H^{(2)}(\varPi)}=\mathrm{tr}\varPi f\varPi\overline h\varPi=\langle\varPi f\varPi, \varPi h\varPi \rangle_{HS}.
\end{equation*}

By definition, we have

\begin{equation*}
\langle f,h\rangle_{H^{(2)}(\varPi)}=\iint\limits_{EE}f(x)|\varPi(x,y)|^2\overline {h(y)}d\mu(x)d\mu(y)=\|\sqrt{f}\varPi\sqrt{h}\|^2.
\end{equation*}
\begin{remark}
The same determinantal measure can be represented by different kernels, but the space $H^{(2)}(\varPi)$ is defined in terms of
the determinantal measure itself, since the measure $|\varPi(x,y)|^2d\mu(x)d\mu(y)$ is uniquely determined by the measure $\mathbb P_{\varPi}$.
\end{remark}
For $r_1 > 0, r_2>0, r_3>0$, let $\mathscr{E}_\Pi(r_1, r_2, r_3)$ be the family of functions ${f \in L_\infty(E, \mu)}$
defined by the formula
\begin{multline*} \mathscr{E}_\Pi(r_1, r_2, r_3) =\\= \{ f \in L_\infty(E, \mu) : \|f\|_{L_\infty} \leq r_1,  \|f \Pi\|_{HS} \leq r_2,
 \left|\det\nolimits_2(1 + \Pi f \Pi)\right|\geq r_3 \}. \end{multline*}

For $\delta > 0$, let $\mathscr{N}_\Pi(\delta)$ be the family of pairs of functions $f \in L_\infty(E, \mu)$ and $h \in L_\infty(E, \mu)$,
defined by the formula
\begin{multline*}
	 \mathscr{N}_\Pi(\delta) = \{ (f, h) : f \in L_\infty(E, \mu), h \in L_\infty(E, \mu),\\
	 \|\Pi f h \Pi\|_{HS} \le \delta, \|\Pi f \Pi h \Pi\|_{HS} \le \delta \}.
\end{multline*}

\begin{theorem}\label{main-thm}
For every $r_1 > 0, r_2>0, r_3>0$, there exist $\delta > 0$ and $A>0$ such that the following holds.
For any locally trace-class projection $\varPi$ and the corresponding determinantal point process $\mathbb P_\varPi$, any $f, h \in \mathscr{E}_\Pi(r_1, r_2, r_3)$ and $(f, h) \in \mathscr{N}_\Pi(\delta)$, we have
\begin{multline*}
	 \left| \mathbb{E}_{\mathbb P_{\varPi}} \widetilde{\Psi}_{1+f} \cdot \widetilde{\Psi}_{1+h} - \mathbb{E}_{\mathbb P_{\varPi}} \widetilde{\Psi}_{1+f} \cdot \mathbb{E}_{\mathbb P_{\varPi}} \widetilde{\Psi}_{1+h} \right| \le\\
	  A \left(  \langle |f|, |h| \rangle_{L_2(\Pi)} + \langle |f|, |h| \rangle_{H^{(2)}(\Pi)} \right).
\end{multline*}
\end{theorem}

\subsection{Correlations of multiplicative functionals of stationary determinantal point processes}
Theorem \ref{main-thm} yields estimates on the rate of mixing for stationary determinantal point processes.
Our convention for the Fourier transform is as follows:
\[
\hat{f}(\lambda) = \int\limits_{\mathbb{R}} e^{-i\lambda t} f(t)\,dt;
\]
\[
f(t) = \frac{1}{2\pi} \int\limits_{\mathbb{R}} \hat{f}(\lambda) e^{i\lambda t}\,d\lambda.
\]
As usual, we denote the inverse Fourier transform of $f$ by $\check{f}$.

We now consider stationary determinantal point processes on $\mathbb{R}$ corresponding to Fourier projections onto bounded Borel sets.
Let $B \subset \mathbb{R}$ be a bounded Borel set of positive measure, and let the projection $\Pi^B$ be defined by the formula
\[
\Pi^B f = (1_B \hat{f})^\vee.
\]
By definition, the kernel of the operator $\Pi^B$ is given by the formula
\[
\Pi^B(x,y) = \check{1}_B(x-y).
\]
In particular, for $B = [-\pi, \pi]$, we obtain the sine process, a determinantal point process with kernel
$$\mathcal S(x,y)=\frac{\sin\pi(x-y)}{\pi(x-y)},$$
the kernel of the projection onto the Paley–Wiener space. Let $g_s$ denote the action of the translation group on the real line, $g_s(t)=t+s$; by a slight abuse of notation,
we use the same symbol $g_s$ for the induced action of the translation group on $\operatorname{Conf}(\mathbb R)$. Recall that
the tail $\sigma$-algebra of a determinantal process associated with an orthogonal projection is trivial (see \cite{Lyons} in the discrete case and \cite{Bufetov-Qiu} and \cite{Osada-Osada} in the general case). In particular, the determinantal process ${\mathbb P_{\Pi^B}}$ also has a trivial tail $\sigma$-algebra. Consequently, the translation group $g_s$ induces a Kolmogorov flow with respect to the determinantal measure ${\mathbb P_{ \Pi^B}}$ and, consequently, has Lebesgue spectrum of countable multiplicity. Our aim is to give explicit regularity estimates for the spectral measures of multiplicative functionals with respect to the determinantal measure ${\mathbb P_{ \Pi^B}}$.

To a bounded Borel set $B$ assign a positive increasing function $\Phi_B$ by setting
\[
\Phi_B(t) = \left( \sup_{r: |r| \geqslant |t|} |\check{1}_B(r)|^2 \right)^{-1}.
\]
Note that positivity of $\Phi_B$ follows from the Heisenberg uncertainty principle (the Fourier transform of the indicator function $1_B(\cdot)$ is not compactly supported, and hence $\Phi_B(t)$ never takes the value $+\infty$). Our key assumption is
\begin{equation}\label{int-phib}
\int\limits_{0}^{+\infty} \frac{dt}{\Phi_B(t)} < +\infty.
\end{equation}

Introduce the class $\mathscr{K}_B$ of bounded Borel functions $f$ on $\mathbb{R}$ such that
\[
\sup_{t \in \mathbb{R}} \Phi_B(t)|f(t)| < +\infty
\]
and set $\|f\|_{\mathscr{K}_B} = \sup\limits_{t \in \mathbb{R}} \Phi_B(t)|f(t)|$.

\begin{corollary} \label{corr-stat}Let $B \subset \mathbb{R}$ be a bounded Borel set satisfying condition \eqref{int-phib}.
For every $C > 0$ and $c>0$, there exists $A > 0$ such that for any Borel functions $f_i\in \mathscr{K}_B$, $i=1,2$, such that
\begin{equation}\label{u1-gen}
    \left|\det\nolimits_2\left(1+\Pi^Bf_i\Pi^B\right)\right|\geq c, \quad \|f_i\|_{\mathscr{K}_B}\leq C,
\end{equation}
the following inequality holds:
$$\sup\limits_{t\in\mathbb{R}}\Phi_B(t/4) |\mathbb{E}_{\mathbb P_{ \Pi^B}}\left({\Psi}_{1+f_1}\right)\circ g_t\cdot {\Psi}_{1+f_2}-\mathbb{E}_{\mathbb P_{ \Pi^B}}{\Psi}_{1+f_1}\mathbb{E}_{\mathbb P_{ \Pi^B}}{\Psi}_{1+f_2}|\leq A.$$
\end{corollary}

In the particular case of the sine process, the expression
$$
(1+t^2)\Phi_{[-\pi, \pi]}(t)
$$
is bounded away from zero and bounded above by definition, and we obtain the following.

\begin{corollary}\label{corr-sine}
For every $C > 0$ and $c>0$, there exists $A > 0$ such that if the Borel functions $f_i:\mathbb R\to\mathbb{C}$, $i=1,2$, satisfy the inequalities

\begin{equation}\label{u1}
    \left|\det\nolimits_2\left(1+{\mathcal S}f_i{\mathcal S}\right)\right|\geq c, \  \sup\limits_{t\in\mathbb{R}}(1+t^2)|f_i(t)|\leq C, \quad i=1,2,
\end{equation}
then the following inequality also holds:
$$\sup\limits_{t\in\mathbb{R}}(1+t^2)|\mathbb{E}_{\mathbb P_{\mathcal S}}\left({\Psi}_{1+f_1}\right)\circ g_t\cdot {\Psi}_{1+f_2}-\mathbb{E}_{\mathbb P_{\mathcal S}}{\Psi}_{1+f_1}\mathbb{E}_{\mathbb P_{\mathcal S}}{\Psi}_{1+f_2}|\leq A.$$
\end{corollary}
Thus, the multiplicative functional corresponding to a function $f$ satisfying
conditions \eqref{u1} has, under the sine process, continuous spectral density with square-integrable derivative; in other words, the spectral density belongs to the Sobolev space $H^1$.

The proofs of Corollaries \ref{corr-stat} and \ref{corr-sine} are given in Section \ref{sec-stat}.

\subsection{Outline of the proof of Theorem \ref{main-thm}}

It will be more convenient for us to work with operator norms rather than determinants and to estimate the multiplicative rather than the additive correlation. We therefore reformulate Theorem \ref{main-thm} as follows.
For $r > 0$, let $\mathscr{C}_\Pi(r)$ be the family of functions $f \in L_\infty(E, \mu)$,
defined by
$$ \mathscr{C}_\Pi(r) = \{ f \in L_\infty(E, \mu) : \|f\|_{L_\infty} + \|f \Pi\|_{HS} + \|(1 + \Pi f \Pi)^{-1}\| < r \}. $$
\begin{proposition}\label{mult-prop}
For every $r > 0$, there exist $\delta > 0$ and $A>0$ such that the following holds.
For any locally trace-class projection $\varPi$ and the corresponding determinantal point process $\mathbb P_\varPi$,
any $f, h \in \mathscr C_\varPi(r)$ and $(f, h) \in \mathscr{N}_\varPi(\delta)$, we have
		\begin{multline}\label{mult-est-form}
		\exp\left( -A(\langle |f|, |h|\rangle_{H^{(2)}(\varPi)}+\langle |f|, |h|\rangle_{L_2(\varPi)}) \right) \le\\
        \Big|\frac{\mathbb{E}_{\mathbb P_{\varPi}} \widetilde{\Psi}_{1+f} \cdot \widetilde{\Psi}_{1+h}}{\mathbb{E}_{\mathbb P_{\varPi}} \widetilde{\Psi}_{1+f} \cdot \mathbb{E}_{\mathbb P_{\varPi}} \widetilde{\Psi}_{1+h}}\Bigg| \\
		\le \exp\left( A(\langle |f|, |h|\rangle_{H^{(2)}(\varPi)}+\langle |f|, |h|\rangle_{L_2(\varPi)}) \right).
	\end{multline}
	\end{proposition}

Theorem \ref{main-thm} immediately follows from Proposition \ref{mult-prop} and the following simple observation.

\begin{proposition}\label{corr-e}For any $r_1 > 0, r_2>0, r_3>0$, there exists $r>0$ such that
$$\mathscr C_\varPi(r)\subset \mathscr E_\varPi(r_1, r_2, r_3).$$ For any $r>0$, there exist $r_1 > 0, r_2>0, r_3>0$ such that
$$  \mathscr E_\varPi(r_1, r_2, r_3)\subset \mathscr C_\varPi(r).$$
\end{proposition}

Thus, under the assumptions of Theorem \ref{main-thm}, there exists $c_1$, depending only on $r$, such that
\[
|\mathbb{E}_{\mathbb P_{\varPi}} \widetilde{\Psi}_{1+f}|^{-1} \geq c_1, \quad |\mathbb{E}_{\mathbb P_{\varPi}} \widetilde{\Psi}_{1+h}|^{-1} \geq c_1.
\]
Estimate \eqref{mult-est-form} now implies the existence of a positive constant $A$ such that
\[
\left| \mathbb{E}_{\mathbb P_{\varPi}} \widetilde{\Psi}_{1+f}\widetilde{\Psi}_{1+h} - \mathbb{E}_{\mathbb P_{\varPi}} \widetilde{\Psi}_{1+f} \mathbb{E}_{\mathbb P_{\varPi}} \widetilde{\Psi}_{1+h}
\right|
\]
\[
\leq A \left( \langle |f|, |h| \rangle_{H^{(2)}(\varPi)} + \langle |f|,|h|\rangle_{L_2(E,\varPi d\mu)} \right).
\]
Theorem \ref{main-thm} follows now from Proposition \ref{mult-prop}. We turn to the proof of Proposition \ref{mult-prop}.

Our first step is to represent the correlation ratio as a Fredholm determinant.

\begin{lemma}\label{corr-fred}
Let $f,h\in HS(\varPi)$ such that the operators $1 + \varPi f \varPi$ and $1 + \varPi h \varPi$ are invertible. Then
\begin{multline}\label{eq-corr-fred}
\frac{\mathbb{E}_{\mathbb P_{\varPi}}\widetilde{\Psi}_{1+f} \cdot\widetilde{\Psi}_{1+h}}{\mathbb{E}_{\mathbb P_{\varPi}}\widetilde{\Psi}_{1+f} \cdot \mathbb{E}_{\mathbb P_{\varPi}}\widetilde{\Psi}_{1+h}} =\\
=\det\bigl(1 + (\varPi fh\varPi - \varPi f\varPi h\varPi) (1 + \varPi h\varPi)^{-1}(1 + \varPi f\varPi)^{-1}\bigr).
\end{multline}
\end{lemma}
The proof of Lemma \ref{corr-fred} is given in Section \ref{sec-fred}.

Set
\begin{equation}\label{def-kfh}
K_\varPi(f, h) = (\varPi fh\varPi - \varPi f\varPi h\varPi) (1 + \varPi h\varPi)^{-1} (1 + \varPi f\varPi)^{-1}.
\end{equation}

Rewrite equation \eqref{eq-corr-fred} in the formula

\begin{equation}\label{fred-kfh}
\frac{\mathbb{E}_{\mathbb P_{\varPi}}\widetilde{\Psi}_{1+f} \cdot\widetilde{\Psi}_{1+h}}{\mathbb{E}_{\mathbb P_{\varPi}}\widetilde{\Psi}_{1+f} \cdot \mathbb{E}_{\mathbb P_{\varPi}}\widetilde{\Psi}_{1+h}}
=\det\bigl(1 + K_\varPi(f, h)\bigr).
\end{equation}

To estimate the Fredholm determinant $\det\bigl(1 + K_\varPi(f, h)\bigr)$, we first estimate the regularized determinant \[
\det\nolimits_2(1+K_\varPi(f, h))=\det\nolimits_2\bigl(1 + (\varPi fh\varPi - \varPi f\varPi h\varPi) (1 + \varPi h\varPi)^{-1}(1 + \varPi f\varPi)^{-1}\bigr),
\] and then estimate the trace $\mathrm{tr}(K_\varPi(f,h))$ separately.

The following lemma gives upper and lower bounds for the regularized Hilbert–Carleman determinant $\det\nolimits_2(1+K_\varPi(f, h))$.

\begin{lemma}\label{det-est-sup}
\begin{enumerate}
\item
For every $r > 0$, there exists $C_{2}(r) > 0$ such that, for all $f, h \in \mathscr{C}_\Pi(r)$, the following inequality holds:
\begin{multline}
	\left|\det\nolimits_2(1 + K_\varPi(f,h) )\right| \le\\ \le
\exp\left(C_{2}(r)\left(\langle |f|^2, |h|^2 \rangle_{L_2(\Pi)} + \langle |f|^2, |h|^2 \rangle_{H^{(2)}(\Pi)}\right)\right).
\end{multline}

\item For every $r > 0$ and $\delta_0>0$, there exists $C_{3}(r, \delta_0) > 0$ such that for all $f, h \in \mathscr{C}_\Pi(r)$ with $(f, h) \in \mathscr{N}_\Pi(\delta_0)$, the following inequality holds:
\begin{multline}
	 \left|\det\nolimits_2(1 + (\Pi f h \Pi - \Pi f \Pi h \Pi) R_h R_f)\right| \ge \\ \ge
 \exp\left(-C_{3}(r)\left(\langle |f|^2, |h|^2 \rangle_{L_2(\Pi)} + \langle |f|^2, |h|^2 \rangle_{H^{(2)}(\Pi)}\right)\right).
\end{multline}
\end{enumerate}
\end{lemma}
Note that the upper bound requires only that $f,h$ belong to $\mathscr{C}_\Pi(r)$, whereas the lower bound relies on the condition $(f,h) \in \mathscr{N}_\Pi(\delta_0)$.
The proof of Lemma \ref{det-est-sup} is given in Section \ref{sec-fred}.

After estimating the regularized Hilbert–Carleman determinant $${\det\nolimits_2(1+K_\varPi(f, h)),}$$ in order to prove Proposition \ref{mult-prop} and the main Theorem \ref{main-thm}, it remains to estimate the trace $\mathrm{tr}(K_\varPi(f,h))$.

\begin{lemma}\label{est-trace}
For any $r > 0$, there exists $C_{4}(r) > 0$, depending only on $r$, such that for all $f, h \in \mathscr{C}_{\varPi}(r)$, we have
		\begin{equation}
		\left| \mathrm{tr} (K_\varPi(f,h))-\langle f, \overline{h}\rangle_{L_2(\varPi)}  \right| \le
        C_{4}(r) \langle |f|,|h|\rangle_{H^{(2)}(\varPi)}.
		\end{equation}
	\end{lemma}

The proof of Lemma \ref{est-trace} is given in Section \ref{sec-trace}.
Lemmas \ref{det-est-sup} and \ref{est-trace} directly imply Proposition \ref{mult-prop}, which in turn implies Theorem \ref{main-thm}.

\section{Proofs of Lemmas \ref{corr-fred} and \ref{det-est-sup}}\label{sec-fred}
\subsection{Proof of Lemma \ref{corr-fred}}
 Lemma \ref{corr-fred} can be proved as follows: a direct computation gives
\begin{multline}\label{eq1}
\frac{\mathbb{E}_{\mathbb P_{\varPi}} \widetilde{\Psi}_{1+f} \cdot \widetilde{\Psi}_{1+h}}{\mathbb{E}_{\mathbb P_{\varPi}} \widetilde{\Psi}_{1+f} \cdot \mathbb{E}_{\mathbb P_{\varPi}} \widetilde{\Psi}_{1+h}} =\frac{\mathbb{E}_{\mathbb P_{\varPi}} \Psi_{1+f} \cdot \Psi_{1+h}}{\mathbb{E}_{\mathbb P_{\varPi}} \Psi_{1+f} \cdot \mathbb{E}_{\mathbb P_{\varPi}} \Psi_{1+h}}=\\
 =\det(1 + (\Pi f h \Pi - \Pi f \Pi h \Pi)(1+\Pi h \Pi)^{-1}(1+\Pi f \Pi)^{-1})
\end{multline}
for $f,h$ such that the operators $\Pi f \Pi$ and $\Pi h \Pi$ are trace class. In the general case, when $f,h$ satisfy $f \Pi \in HS$ and $h \Pi \in HS$, the formula follows by continuity.

Equality \eqref{eq1} can also be verified by direct calculation as follows. By definition, we have
$$ \widetilde{\Psi}_{1+h} \cdot \widetilde{\Psi}_{1+f} = \widetilde{\Psi}_{(1+h)(1+f)} \times \exp(\mathrm{tr} \varPi hf\varPi). $$

Recall also that for any Hilbert–Schmidt operators $A,B$, the following identity holds:
$$ \det\nolimits_2((1+A)(1+B)) = \det\nolimits_2(1+A) \cdot \det\nolimits_2(1+B) \cdot \exp(- \mathrm{tr}(AB)).$$

We can therefore write
\begin{multline*}
	\frac{\mathbb{E}_{\mathbb P_{\varPi}}\widetilde{\Psi}_{1+f} \cdot \widetilde{\Psi}_{1+h}}{\mathbb{E}_{\mathbb P_{\varPi}}\widetilde{\Psi}_{1+f} \cdot {\mathbb{E}_{\mathbb P_{\varPi}}\widetilde{\Psi}_{1+h}}} = \\
 = \frac{\det\nolimits_2(1+\varPi(f+h)\varPi + \varPi fh\varPi)}{\det\nolimits_2(1+\varPi(f+h)\varPi + \varPi f\varPi h\varPi)} \times \exp \mathrm{tr}(\varPi fh\varPi - \varPi f\varPi h\varPi).
\end{multline*}

Furthermore, for any Hilbert–Schmidt operator $A$ such that $1+A$ is invertible, and any Hilbert–Schmidt operator $\Delta$, we have the identity
$$ \frac{\det\nolimits_2(1+A+\Delta)}{\det\nolimits_2(1+A)} = \det\nolimits_2\left(1+\Delta(1+A)^{-1}\right) \cdot \exp(- \mathrm{tr}\left(\Delta A(1+A)^{-1}\right)). $$

Applying this identity, we obtain
\begin{multline*} \frac{\mathbb{E}_{\mathbb P_{\varPi}}\widetilde{\Psi}_{1+h} \cdot \widetilde{\Psi}_{1+f}}{\mathbb{E}_{\mathbb P_{\varPi}}\widetilde{\Psi}_{1+h}\cdot\mathbb{E}_{\mathbb P_{\varPi}}\widetilde{\Psi}_{1+f}} = \\
= \det\nolimits_2\left(1 + (\varPi hf\varPi -\varPi h\varPi f\varPi)(1+\varPi h\varPi)^{-1}(1+\varPi f\varPi)^{-1}\right) \times \\
 \times \exp\left(\mathrm{tr}(\varPi hf\varPi -\varPi h\varPi f\varPi)(1+\varPi h\varPi)^{-1}(1+\varPi f\varPi)^{-1}\right) \\
=\det\bigl(1 + (\varPi fh\varPi - \varPi f\varPi h\varPi) (1 + \varPi h\varPi)^{-1}(1 + \varPi f\varPi)^{-1}\bigr).
\end{multline*}
 Lemma \ref{corr-fred} is proved.
$\blacksquare$
\begin{remark}
Recall here that multiplying a determinantal measure by a regularized multiplicative functional again yields a determinantal point process; see \cite{Bufetov12}.
\end{remark}
\subsection{Proof of Lemma \ref{det-est-sup}}

Set
\begin{equation}
 R_h =(1 + \varPi h\varPi)^{-1},
 R_f =(1 + \varPi f\varPi)^{-1}.
\end{equation}

We have
\begin{multline*}
 R_h = 1 - \Pi h \Pi R_h = 1 - R_h \Pi h\Pi, \
 R_f = 1 - \Pi f \Pi R_f = 1 - R_f \Pi f\Pi.
\end{multline*}

By assumption, we have
\begin{equation}\label{r-cond}
\|R_f\|\leq r, \|R_h\|\leq r.
\end{equation}

Our aim is to estimate \[
\det\nolimits_2(1+K_\varPi(f, h))=\det\nolimits_2\bigl(1 + (\varPi fh\varPi - \varPi f\varPi h\varPi) (1 + \varPi h\varPi)^{-1}(1 + \varPi f\varPi)^{-1}\bigr).
\]
To obtain the upper bound, it suffices to bound the Hilbert–Schmidt norm of the operator $(\varPi fh\varPi - \varPi f\varPi h\varPi)$ from above and apply the classical Carleman inequality. Indeed, the upper bound follows directly from Carleman's inequality (cf. \eqref{carleman-ineq}), the identity
$$ \|\Pi f h \Pi\|_{HS}^2 = \langle |f|^2, |h|^2 \rangle_{L_2(\Pi)} $$
and the inequality $\|\Pi f \Pi h \Pi\|_{HS}^2 \le \langle |f|^2, |h|^2 \rangle_{H^{(2)}(\Pi)}$.
To obtain the lower bound, apply Lemma \ref{inf-est-hilb-det} to $A=K_\varPi(f,h)$. In view of condition \eqref{r-cond}, for every $\delta_1>0$ there exists $\delta_0>0$ such that $(f,h)\in\mathscr{N}_\Pi(\delta_0)$ implies $\|K_\varPi(f,h)\|_{HS}<\delta_1$ and, a fortiori, $\|K_\varPi(f,h)\|_{HS}<\delta_1$. The second part
of Lemma \ref{det-est-sup} now follows directly from Lemma \ref{inf-est-hilb-det}.

\section{Estimating the trace in terms of Hilbert–Schmidt norms and proof of Lemma \ref{est-trace}}\label{sec-trace}
	\begin{proposition}\label{trace-ineq-trace}
For every $r > 0$, there exists $C_{4}(r) > 0$, depending only on $r$, such that for all $f, h \in \mathscr{C}_{\varPi}(r)$, we have
		\begin{multline*}
		\left| \mathrm{tr} ((\Pi f h \Pi - \Pi f \Pi h \Pi) R_h R_f -  \Pi f h \Pi + \Pi f \Pi h \Pi + \Pi f h \Pi (f+h) \Pi) \right| \le \\
        C_{4}(r) \| \sqrt{f} \Pi \sqrt{h} \|_{HS}^2=C_{4}(r)\langle |f|,|h|\rangle_{H^{(2)}(\varPi)}.
		\end{multline*}
	\end{proposition}
By definition,
	  \begin{multline*}
		\mathrm{tr}(\Pi f h \Pi)=\langle f, \overline{h}\rangle_{L_2(\varPi)},\\
		\mathrm{tr}(\Pi f \Pi h \Pi)=\langle f, \overline{h} \rangle_{H^{(2)}(\varPi)}, \\
	\mathrm{tr} (\Pi f h \Pi (f+h) \Pi)=\langle fh, \overline{f}+\overline{h}\rangle_{H^{(2)}(\varPi)}.
	  \end{multline*}
It is clear that
\[
\bigl|\langle fh,\,\overline{f}+\overline{h}\rangle_{H^{(2)}(\Pi)}\bigr|
\leqslant \bigl(|f|_{L_\infty} + |h|_{L_\infty}\bigr)\,
\bigl\langle |f|,\, |h|\bigr\rangle_{H^{(2)}(\Pi)}.
\]

Thus, Proposition \ref{trace-ineq-trace} yields

\begin{corollary}\label{corr-trace-fin}
For every $r>0$, there exists $C_{5}(r)>0$, depending only on $r$, such that for all $f,h \in \mathcal{C}_\Pi(r)$, we have
\[
\bigl|\operatorname{tr}\bigl((\Pi f h \Pi - \Pi f \Pi h \Pi)\,R_h R_f\bigr) - \langle f, \overline{h}\rangle_{L_2(\Pi)}\bigr|
\leqslant C_{5}(r)\,\bigl\langle |f|,\, |h|\bigr\rangle_{H^{(2)}(\Pi)}.
\]
\end{corollary}

Corollary \ref{corr-trace-fin} follows directly from Proposition \ref{trace-ineq-trace} and it in turn implies Lemma \ref{est-trace}. We proceed to the proof of Proposition \ref{trace-ineq-trace}.
The proof of Proposition \ref{trace-ineq-trace} is based on the following simple observation: if $Q: L_2(E, \mu) \to L_2(E, \mu)$ is a bounded operator, then we have
\begin{multline}
|\operatorname{tr} \Pi f \Pi h Q h| \leq \|Q\| \cdot \|\sqrt{f}\Pi h\|_{HS}^2, \\
|\operatorname{tr} \Pi h \Pi f Q f| \leq \|Q\| \cdot \|\sqrt{h}\Pi f\|_{HS}^2, \\
|\operatorname{tr} \Pi f h \Pi h Q h| \leq \|Q\| \cdot \|h\|_{L_\infty} \cdot \|\sqrt{f}\Pi h\|_{HS}^2, \\
|\operatorname{tr} \Pi f h \Pi h Q f| \leq \|Q\| \cdot \|f\|_{L_\infty} \cdot \|\sqrt{h}\Pi f\|_{HS}^2.
\end{multline} It is clear that
\[
\|\sqrt{f}\,\Pi h\|_{HS}^2 \leqslant \|h\|_{L_\infty} \cdot \|\sqrt{f}\,\Pi\sqrt{h}\|_{HS}^2,
\]
and similarly for the norm $\|\sqrt{h}\,\Pi f\|_{HS}^2$.

\smallskip 

Write
	\begin{multline*}
		(\Pi f h \Pi - \Pi f \Pi h \Pi) R_h R_f =
		 \Pi f h \Pi - \Pi f \Pi h \Pi - \\
		 - \Pi f h \Pi h \Pi R_h - \Pi f h \Pi f \Pi R_f +
		  \Pi f \Pi h \Pi h \Pi R_h + \Pi f \Pi h \Pi f \Pi R_f + \\
		 + (\Pi f h \Pi - \Pi f \Pi h \Pi) \Pi h \Pi R_h R_f \Pi f \Pi.
	\end{multline*}
	
It is clear that
	\begin{equation*}
		|\mathrm{tr} (\Pi f h \Pi - \Pi f \Pi h \Pi) \Pi h \Pi R_h R_f \Pi f \Pi | \le C_{6}(r) \| \sqrt{f} \Pi \sqrt{h} \|_{HS}^2.
	\end{equation*}
It follows from the above estimates that we have
	\begin{equation*}
		|\mathrm{tr} \Pi f \Pi h \Pi h \Pi R_h | = | \mathrm{tr} \Pi f \Pi h \Pi R_h \Pi h \Pi | \le C_{7}(r) \| \sqrt{f} \Pi \sqrt{h} \|_{HS}^2,
	\end{equation*}
	\begin{equation*}
		|\mathrm{tr} \Pi f h \Pi f \Pi R_f | = | \mathrm{tr} \Pi f h \Pi R_f \Pi f \Pi | \le C_{8}(r) \| \sqrt{f} \Pi \sqrt{h} \|_{HS}^2
	\end{equation*}

Next, we have
\begin{equation*}
	\Pi f h \Pi R_h = \Pi f h \Pi h \Pi - \Pi f h \Pi h \Pi R_h \Pi h \Pi.
\end{equation*}

It follows from the above estimates that we have
\[
|\operatorname{tr}\, \Pi f h \Pi h \Pi\, R_h R_f\, \Pi f \Pi|
\leqslant C_{9}(r)\,\|\sqrt{f}\,\Pi\sqrt{h}\|_{HS}^2;
\]
\begin{multline*}
|\operatorname{tr}\, \Pi f \Pi h \Pi h \Pi\, R_h R_f\, \Pi f \Pi|
= |\operatorname{tr}\, \Pi f \Pi h \Pi\, R_h\, \Pi h \Pi f \Pi\, R_f|\leqslant \\
\leqslant C_{10}(r)\,\|f \Pi h\|_{HS}^2
\leqslant C_{10}(r)\,\|f\|_{L_\infty}^2 \cdot \|h\|_{L_\infty}^2 \cdot \|\sqrt{f}\,\Pi\sqrt{h}\|_{HS}^2,
\end{multline*}
which in turn implies the inequality

\begin{equation*}
	| \mathrm{tr} \Pi f h \Pi h \Pi R_h \Pi h \Pi | \le C_{11}(r) \| \sqrt{f} \Pi \sqrt{h} \|_{HS}^2.
\end{equation*}

Next, we have the identity
\begin{equation*}
	\Pi f h \Pi f \Pi R_f = \Pi f h \Pi f \Pi - \Pi f h \Pi f \Pi R_f \Pi f \Pi,
\end{equation*}

and we clearly have the inequality
\begin{equation*}
	| \mathrm{tr} \Pi f h \Pi f \Pi R_f \Pi f \Pi | \le C_{12}(r) \| \sqrt{f} \Pi \sqrt{h} \|_{HS}^2.
\end{equation*}
	
Proposition \ref{trace-ineq-trace} is proved.
$\blacksquare$

\section{Spectral density of multiplicative functionals of stationary determinantal point processes}\label{sec-stat}
\subsection{Proofs of Corollaries \ref{corr-stat} and \ref{corr-sine}}

We will need a simple lemma on the convolution of rapidly decaying functions.

Let $\Phi: [0, +\infty) \to (0, +\infty)$ be a positive increasing function satisfying
\[
\int\limits_{0}^{+\infty} \frac{dt}{\Phi(t)} < +\infty.
\]
Extend $\Phi$ to an even function by setting $\Phi(-t) = \Phi(t)$ for $t>0$.

\begin{lemma}\label{conv-ineq} Let $f,h: \mathbb{R} \to \mathbb{R}$ be Borel functions such that
\[
\sup_{t \in \mathbb{R}} \Phi(t)(|f(t)| + |h(t)|) = A < +\infty.
\]
We then have the inequality
\[
\sup_{R \in \mathbb{R}} \Phi\left(\frac{R}{2}\right) \cdot |f * h(R)| \leqslant A \left(\|f\|_{L_1} + \|h\|_{L_1}\right).
\]
\end{lemma}

 Proof of Lemma \ref{conv-ineq}. Let $R > 0$. We write
\begin{align*}
f * h(R) &= \int\limits_{-\infty}^{+\infty} f(R-u)h(u)\,du = \\
&= \int\limits_{-\infty}^{R/2} f(R-u)h(u)\,du + \int\limits_{R/2}^{+\infty} f(R-u)h(u)\,du.
\end{align*}
The absolute value of the first integral is at most $A\|h\|_{L_1}/\Phi(R/2)$, while the absolute value of the second is at most $A\|f\|_{L_1}/\Phi(R/2)$. The proof of Lemma \ref{conv-ineq} is now complete.
$\blacksquare$

Next, for $\alpha > 0$, let $\mathcal{K}_{B,\alpha}$ be the set of Borel functions $f$ on $\mathbb{R}$ such that $\sup_{t \in \mathbb{R}} \Phi_B(\alpha t)|f(t)| < +\infty$, and define
\[
\|f\|_{\mathcal{K}_{B,\alpha}} = \sup_{t \in \mathbb{R}} \Phi_B(\alpha t)|f(t)|.
\]

Recall that we make the key condition (\ref{int-phib}), as above. For any $f,h\in\mathcal{K}_B$, we now have $f*h\in\mathcal{K}_{B,1/2}$, and
\[
\|f * h\|_{\mathcal{K}_{B,1/2}} \leqslant \bigl(\|f\|_{\mathcal{K}_B} + \|h\|_{\mathcal{K}_B}\bigr) \cdot \left(\int_{\mathbb{R}} \frac{dt}{\Phi_B(t)}\right).
\]

The following simple proposition, which estimates the inner product $${\langle f,h(\cdot+R)\rangle_{H^{(2)}(\Pi_B)}}$$ in terms of $\Phi_B$, plays a key r\^ole in deriving Corollaries \ref{corr-stat} and \ref{corr-sine} from Theorem \ref{main-thm}.

\begin{proposition}\label{h2norm-est}If a bounded Borel set $B$ satisfies assumption (\ref{int-phib}),
then for all $f,h\in\mathcal{K}_B$ we have
\[
 \Phi_B\!\left(\frac{R}{4}\right)\cdot \bigl|\langle f, h(\cdot+R)\rangle_{H^{(2)}(\Pi^B)}\bigr| \leqslant 4\bigl(\|f\|_{\mathcal{K}_B} + \|h\|_{\mathcal{K}_B}\bigr)\bigl(1 + \|\Phi_B^{-1}\|_{L_1}\bigr).
\]
\end{proposition}
 Proof of Proposition \ref{h2norm-est}. By definition,
\[
\langle f, h(\cdot+R) \rangle_{H^{(2)}(\Pi^B)} = \int\limits_{-\infty}^{+\infty}\!\int\limits_{-\infty}^{+\infty} f(u)\,h(u+R+k)\,|\Pi^B(k)|^2\,du\,dk.
\]
Define $\tilde{f}$ by the formula $\tilde{f}(k)=f(-k)$. Set $g=\tilde{f}*h$. By Lemma \ref{conv-ineq}, $g\in\mathcal{K}_{B,1/2}$ and we have
\[
\|g\|_{\mathcal{K}_{B, 1/2}} \leqslant \bigl(\|f\|_{\mathcal{K}_B} + \|h\|_{\mathcal{K}_B}\bigr) \cdot \|\Phi_B^{-1}\|_{L_1}.
\]
Write $$\langle f,h(\cdot+R)\rangle=\int_{\mathbb{R}}g(R+k)\,|\Pi^B(k)|^2\,dk.$$ Applying Lemma \ref{conv-ineq} once again proves the proposition.

Finally, note that if $f,h\in\mathcal{K}_B$, then for every $\delta>0$, the condition $$(f,h(\cdot+R))\in\mathcal{N}_{\Pi^B}(\delta)$$ automatically holds for all sufficiently large $R$.

For the correlation coefficient of multiplicative functionals from  Theorem \ref{main-thm}, we conclude that for any bounded Borel set $B\subset\mathbb{R}$ satisfying condition \eqref{int-phib}, and for
any $r>0$ and $\delta>0$, there exists $A>0$ such that, if $$f,h\in\mathscr{K}_B\bigcap\mathscr{C}_{\Pi^B}(r)$$ and $\|f\|_{\mathscr{K}_B},\|h\|_{\mathscr{K}_B}<D$ then

\begin{equation}\label{C9}
\Phi\!\left(\frac{R}{4}\right)\cdot \bigl|\mathbb E_{P_{\Pi^B}} \widetilde{\Psi}_{1+f} \cdot \widetilde{\Psi}_{1+h(\cdot+R)} -\mathbb E_{P_{\Pi^B}} \widetilde{\Psi}_{1+f} \cdot\mathbb E_{P_{\Pi^B}} \widetilde{\Psi}_{1+h} \bigr| \leqslant A.
\end{equation}

This proves Corollaries \ref{corr-stat} and \ref{corr-sine}.

\subsection{The case of additive functionals}
We compare our estimates for the decay of correlation of multiplicative functionals with the formulas for the spectral densities of additive functionals obtained by A. B. Soshnikov (cf. (3.11) in \cite{Soshnikov-2}). To a bounded Borel function $f$ with compact support on $E$ assign the additive functional $S_f$ on $\mathrm{Conf}(E)$ by the formula \[
S_f(X) = \sum_{x \in X} f(x).
\] We have $$\mathbb{E}S_f = \int_E f(x)\Pi(x,x)d\mu(x).$$ Set $\widetilde{S_f}=S_f-\mathbb{E}S_f$.
A direct calculation gives\[
\mathbb{E}_{P_{\Pi}}|\widetilde{S_f}|^2 = \frac{1}{2} \int\limits_E  \int\limits_E |f(x)-f(y)|^2 |\Pi(x,y)|^2 d\mu(x)d\mu(y).
\] Define the Hilbert space $\dot{H}_{1/2}(\Pi)$ as the completion of the family of all bounded Borel functions with compact support on $E$ with respect to the Euclidean norm $\|f\|_{\dot{H}_{1/2}(\Pi)}$, defined by
\[
\|f\|_{\dot{H}_{1/2}(\Pi)}^2 = \frac{1}{2} \int\limits_E  \int\limits_E |f(x)-f(y)|^2 |\Pi(x,y)|^2 d\mu(x)d\mu(y).
\]

The map $f\to\widetilde{S_f}$ extends to an isometric embedding of the Hilbert space $\dot{H}_{1/2}(\Pi)$ into $L_2(\mathrm{Conf}(E),P_{\Pi})$.

We now return to the special case of a stationary determinantal point process on $\mathbb{R}$ with the kernel
\[
\Pi^B(x,y) = \check{\mathbf{1}}_B(x-y).
\]
A direct calculation gives
\begin{multline}\label{comp-spec-add}
\mathbb{E}_{P_{\Pi^B}}\widetilde{S_f} \cdot \overline{\widetilde{S_{f(\cdot+R)}}} =\\= \frac{1}{2} \int\limits_{\mathbb{R}} \int\limits_{\mathbb{R}} (f(x)-f(y))(\overline{f(x+R)-f(y+R)}) |\Pi^B(x,y)|^2 dx\,dy=  \\
= \frac{1}{2} \int\limits_{\mathbb{R}} \int\limits_{\mathbb{R}}  (f(x+h)-f(x))(\overline{f(x+h+R)-f(x+R)}) |\check{\mathbf{1}}_B(h)|^2 dxdh =\\
= \frac{1}{4\pi} \int\limits_{\mathbb{R}} \int\limits_{\mathbb{R}}  |\hat{f}(\lambda)|^2 |e^{i\lambda h}-1|^2 |\check{\mathbf{1}}_B(h)|^2 e^{-i\lambda R} d\lambda dh =\\
= \frac{1}{2\pi} \int\limits_{\mathbb{R}}  |\hat{f}(\lambda)|^2 (\widehat{|\check{\mathbf{1}}_B|^2}(0) - \widehat{|\check{\mathbf{1}}_B|^2}(\lambda)) e^{-i\lambda R} d\lambda.
\end{multline}
In the last step, we used the fact that the function $|\check{\mathbf{1}}_B|^2$ is even.
Set
\[
Q_B = \widehat{|\check{\mathbf{1}}_B|^2}.
\]
Let $\rho^{B,f}$ stand for the spectral density of the regularized additive functional $\widetilde{S_f}$ with respect to the determinantal measure $P_{\Pi^B}$.

Following A. B. Soshnikov, observe that equation \eqref{comp-spec-add} implies the identity
\[
\rho^{B,f}(\lambda) = |\hat{f}(\lambda)|^2 (Q_B(0) - Q_B(\lambda)).
\]
Using Lemma \ref{conv-ineq}, one can directly verify that if $f\in\mathcal{K}_B$, then we have
\[
\sup_{R \in \mathbb{R}} \Phi_B(R/4) |\mathbb{E}_{P_{\Pi^B}}\widetilde{S_f} \cdot \overline{\widetilde{S_{f(\cdot+R)}}}| < +\infty.
\]

Therefore, Theorem \ref{main-thm} gives the same rate of correlation decay for multiplicative functionals as Soshnikov's spectral formula does for additive functionals.

\textbf{Question.} What lower bounds can be obtained for the correlation ratio of multiplicative functionals of determinantal point processes?

\section{Bounds for Fredholm and Hilbert–Carleman determinants}
For completeness, this section provides proofs of the elementary estimates for Fredholm and Hilbert–Carleman determinants used in the paper. We first recall Carleman's inequality for the regularized Hilbert–Carleman determinant; see, for example, Gohberg, Goldberg, and Krupnik \cite{gohberg}, Theorem 7.4, formula (7.10), p. 69: for every Hilbert-Schmidt operator $A$,

\begin{equation}\label{carleman-ineq}
|\det\nolimits_2(1+A)|  \leq \exp\bigl(\frac 12\|A\|_{HS}^2\bigr).
\end{equation}

We now obtain a simple bound for the Hilbert–Schmidt norm of the inverse operator in terms of the Hilbert–Schmidt norm of the original operator and its Hilbert–Carleman determinant.
\begin{proposition}\label{inv-hs-est}
Let $A$ be a Hilbert–Schmidt operator such that $1+A$ is invertible. Then
$$ \left\| \frac{A}{1+A} \right\|_{HS} \le \frac{\exp((\|A\|_{HS}^2 + 1))}{|\det\nolimits_2(1+A)|}. $$
\end{proposition}

The proof is a straightforward adaptation of the proof of Corollary 5.2 in Simon \cite{Simon}.

Set
$$B = \frac{\left( \frac{A}{(1+A)} \right)^*}{\|\frac{A}{(1+A)}\|_{HS}}.$$
Clearly, $\mathrm{tr} BA(1+A)^{-1} = \left\| \frac{A}{1+A} \right\|_{HS}$.

Consider the function
\begin{multline*}
	F(\lambda) = \det\nolimits_2(1 + A + \lambda B) =\\
	\det\nolimits_2(1+A) \cdot \det\nolimits_2(1 + \lambda B (1+A)^{-1}) \cdot \exp(-\lambda\|\frac{A}{1+A}\|_{HS}).
\end{multline*}
We have
$$ \frac{d F(\lambda)}{d \lambda} \bigg|_{\lambda=0} = -\det\nolimits_2(1+A) \left\| \frac{A}{1+A} \right\|_{HS}. $$
By Cauchy's estimate for the derivative of an entire function,
$$ \left| \frac{d F(\lambda)}{d \lambda} \bigg|_{\lambda=0} \right| \le \exp((\|A\|_{HS}^2 + 1)). $$
Proposition \ref{inv-hs-est} is proved.
$\blacksquare$

Next, let $\mathscr J_1$ be the ideal of trace-class operators acting on the Hilbert space $L_2(E,\mu)$:
\begin{equation*}
\mathscr J_1=\{A:\mathrm{tr}\sqrt{AA^{\star}}<+\infty\}.
\end{equation*}
As usual, write $\|A\|_{\mathscr J_1}=\mathrm{tr}\sqrt{AA^{\star}}$.

Recall that if $A\in\mathscr J_1$,
and the operator $1+A$ is invertible, then we have
\begin{multline}\label{exp-est}
\left|\det \left((1+A)^{-1}\right)\right|
=\left|\det \left(1 - A(1+A)^{-1}\right)\right|\\ \leq  \exp \operatorname{tr} \left(|A(1+A)^{-1}|\right)
\leq \exp \left( \|A\|_{\mathscr J_1} \|(1+A)^{-1}\| \right).
\end{multline}

Similarly, if $A$ is a Hilbert–Schmidt operator and $1+A$ is invertible, then we have
\begin{multline}
\left|\det\nolimits_2 \left((1+A)^{-1}\right)\right| = \left|\det\nolimits_2 \left( 1 - A(1+A)^{-1} \right)\right| \leq \\
\leq\exp \left( \|A(1+A)^{-1}\|_{HS}^2/2 \right)
\leq \exp \left( \|A\|_{HS}^2 \cdot \|(1+A)^{-1}\|^2/2 \right).
\end{multline}
We also have
\begin{equation}
\det\nolimits_2 \left((1+A)^{-1}\right)\det\nolimits_2 \left((1+A)\right)=\exp(-{\mathrm {tr}} A^2(1+A)^{-1}).
\end{equation}
Proposition \ref{inv-hs-est} now yields the following corollary.
\begin{corollary}\label{hs-inverse-est} For every
$R>0$, there exist positive constants $c_{13}$ and $C_{13}$ depending only on $R$ such that, for every Hilbert–Schmidt operator $A$ whose Hilbert–Schmidt norm does not exceed $R$ and such that the operator $1+A$ is invertible, the following inequalities hold:
\begin{equation}
c_{13}\left(\|(1+A)^{-1}\|_{ HS}\right)^{-1}\leq \left|\det\nolimits_2 \left( 1+A \right)\right| \leq C_{13}\left(\|(1+A)^{-1}\|_{ HS}\right)^{-1};
\end{equation}

\begin{equation}
c_{13}\left|\det\nolimits_2 \left( 1+A \right)\right|^{-1}\leq \|(1+A)^{-1}\|_{ HS} \leq C_{13}\left|\det\nolimits_2 \left( 1+A \right)\right|^{-1}.
\end{equation}

\end{corollary}

If we additionally assume that our Hilbert–Schmidt operator is a contraction, then the regularized Hilbert–Carleman determinant admits a lower bound in terms of the Hilbert–Schmidt norm of the operator:

\begin{lemma}\label{inf-est-hilb-det}
Let $A$ be a Hilbert–Schmidt operator such that $\|A\| < 1$.
Then
		$$ \exp\left( -\frac{\|A\|_{HS}^2}{2(1-\|A\|)^2}-\frac{\|A\|_{HS}^2}{2(1-\|A\|)} \right) \le |\det\nolimits_2 (1+A)|. $$
	\end{lemma}

To prove this, write the identity
\[
\det\nolimits_2(1+A) \cdot \det\nolimits_2\!\bigl((1+A)^{-1}\bigr) = \exp\!\Bigl(\operatorname{tr}\bigl(A^2(1+A)^{-1}\bigr)\Bigr),
\]
which follows from the general identity
\[
\det\nolimits_2(1+A) \cdot \det\nolimits_2(1+B) = \det\nolimits_2\!\bigl((1+A)(1+B)\bigr) \cdot \exp\!\bigl(\operatorname{tr}(AB)\bigr)
\]
by substituting $B=-A(1+A)^{-1}$ (so that $1+B=(1+A)^{-1}$).

Under our assumptions, we have
$$\|(1+A)^{-1}\| = \|1 - A + A^2 - \dots\| \le (1 - \|A\|)^{-1}$$
and
$$|\operatorname{tr} A^2(1+A)^{-1}| \le \frac{\|A\|_{HS}^2}{1 - \|A\|}.$$

Furthermore, we have \begin{multline}|\operatorname{det}_2((1+A)^{-1})| = |\operatorname{det}_2(1 - A(1+A)^{-1})| \le \\ \le
\exp(\|\frac{1}{2}A(1+A)^{-1}\|_{HS}^2) \le \exp\left(\frac{\|A\|_{HS}^2}{2(1-\|A\|)^2}\right).\end{multline}

It follows that we have
\begin{multline}
|\operatorname{det}_2(1+A)| = \left|\frac{\exp(-\operatorname{tr} \frac{1}{2}A^2(1+A)^{-1})}{\operatorname{det}_2((1+A)^{-1})}\right| \ge \\ \ge
\exp\!\left(-\frac{\|A\|_{HS}^2}{2(1-\|A\|)^2} - \frac{\|A\|_{HS}^2}{2(1-\|A\|)}\right),
\end{multline}
as required.
    $\blacksquare$

Note added in proof. After this paper was completed, I became aware of the beautiful paper \cite{dumaz} in which speed of mixing is investigated for the sine-beta process.

\textbf{Acknowledgements.} I am deeply grateful to A.S. Arakcheev, M.A. Borovikov and S.M. Gorbunov for fruitful discussions. I am deeply grateful to the anonymous referees whose detailed comments greatly helped me to improve the presentation.

The author is a winner of the ``Leader'' competition of the BASIS Foundation for the Advancement of Mathematics and Theoretical Physics, and would like to thank the Jury and sponsors.

Part of this paper was written while I was a guest at the Sirius Mathematics Centre, and I am deeply grateful to the Centre for its warm hospitality.

\end{document}